\documentclass[11pt,a4paper]{article}

\usepackage{epic,eepic,epsfig,amsmath,amssymb,amsthm}
\usepackage{t1enc}

\usepackage{color}

\usepackage{bbm}

\def\cB{{\mathcal B}}

\def\cF{{\mathcal F}}

\def\cR{{\mathcal R}}

\def\longformule#1#2{
\displaylines{ \qquad{#1} \hfill\cr \hfill {#2} \qquad\cr } }

\def\uapp{u_{{\rm app}}}

\def\virgp{\raise 2pt\hbox{,}}

\renewcommand{\geq}{\geqslant}
\renewcommand{\leq}{\leqslant}

\def\R{{\mathbb R}}

\def\Lam{\Lambda}

\def\virgp{\raise 2pt\hbox{,}}
\def\cdotpv{\raise 2pt\hbox{;}}

\def\1{\mathbbm{1}}

\newtheorem{theorem}{Theorem}[section]

\newtheorem{proposition}[theorem]{Proposition}
\newtheorem{lemma}[theorem]{Lemme}

\newtheorem{pte}[theorem]{Property}

\theoremstyle{remark}
\newtheorem{remark}{Remark}[section]

\theoremstyle{definition}
\newtheorem{definition}{Définition}[section]

\theoremstyle{definition}
\newtheorem{example}{Exemple}[section]

\theoremstyle{definition}

\begin{document}

\title{A few notes on Lorentz spaces}

\author{Claire David}

\maketitle
\centerline{Sorbonne Universit\'e}

\centerline{CNRS, UMR 7598, Laboratoire Jacques-Louis Lions, 4, place Jussieu 75005, Paris, France}



\maketitle

\section{Introduction}

\vskip 1cm
  
In the sequel, we recall and comment some classical results on the non-increasing rearrangement and Lorentz spaces. There are papers in the existing literature that seemed to have been bypassed as regards its contractive property in~$L^p$ spaces. Also, we provide detailed proofs and a few properties that does not seem to arise in the existing literature.  
\vskip 1cm
  
\section{Framework of the study}
\label{Lorentz}

\noindent In this section, we have collected properties on the decreasing rearrangement, and on Lorentz spaces~$ L^{p,q}(\R^d)$,
\mbox{$1 \leq p < + \infty$}, \mbox{$1 \leq q < + \infty$}, $d \geq 1$. We have tried to give the detailed proofs.

\vskip 1cm


\begin{definition}\textbf{Distribution function}\\


\noindent Let us denote by~$(E,\mu)$ a measured space. If one denotes by~$f$ a real-valued, measurable function, $\mu-$finite a.e., we introduce the positive-valued \emph{distribution function} $\mu_f$, defined on $\R^+$, such that, for any positive number~$\lambda$:

$$\mu_f(\lambda)=\mu \left ( \left \lbrace x \, \big |\, |f(x)| > \lambda\right \rbrace \right)$$

\noindent which can also be written as:

$$  \mu_f(\lambda)=\displaystyle \inf_{s >0}  \,\left \lbrace s  \big |\,\mu\left ( \left \lbrace x \, \big |\, |f(x)| > s\right \rbrace \right) \leq  \lambda \right \rbrace$$

\noindent For any strictly positive number $\sigma $, we set:

$$\begin{array}{ccc}  m(\sigma , f  )&= & \text{mes}\left ( \left \lbrace x \, \big |\, |f ( x)|  > \sigma\right \rbrace \right)
\end{array}$$

\end{definition}

\vskip 1cm

\begin{pte}\textbf{Properties of the distribution function}\\

\noindent The distribution function $\mu_f(\lambda)$ is non-increasing, and right-continuous on~$[0,+ \infty[$.\\

\end{pte}

\vskip 1cm

\begin{proof}

Let us consider a positive number~$\lambda_0$. As in~\cite{BennettSharpley}, we set, for any real positive number~$\lambda$:

$$E(\lambda)=\left \lbrace x \,\in\,\R^d \, \bigg |\,  |f(x) | >\lambda \right \rbrace $$

\noindent As~$\lambda$ increases, the sets~$E(\lambda)$ decrease; moreover:

$$E(\lambda_0)=\underset {\lambda > \lambda_0}{\bigcup} E(\lambda)=\underset{n=1}{\overset{+\infty}{ \bigcup   }}\, E\left (\lambda_0+\displaystyle \frac{1}{n} \right)$$

\noindent In order to prove the right-continuity, one requires to show that:

$$\displaystyle \lim_{\lambda \to 0^+} \mu_f\left (\lambda_0+\lambda \right)=\mu_f\left (\lambda_0 \right)$$

\noindent but since the mapping~\mbox{$ \lambda \mapsto \mu_f\left (\lambda  \right)$} is non-increasing, the monotone convergence theorem yields the right-continuity:

$$\displaystyle \lim_{n \to + \infty} \mu_f\left (\lambda_0+\displaystyle \frac{1}{n} \right)=
\displaystyle \lim_{n \to + \infty} \mu \left (E\left (\lambda_0+\displaystyle \frac{1}{n} \right)\right)=
\mu \left (E\left (\lambda_0  \right)\right)=\mu_f\left (\lambda_0 \right)$$

\end{proof}
\vskip 1cm

\begin{definition}\textbf{Equimeasurable functions}\\

\noindent Let us denote by~$(E,\mu)$ and~$(F,\nu)$ two measured spaces,~$f$ a real-valued function defined on $E$,~$g$ a real-valued function defined on~$F$. The functions~$f$ and~$g$ are said to be equimeasurable if:

  $$ \forall\,t>0\, : \quad   \mu\, \left \lbrace   x\, \in \,E\, : \quad |f(x)| > t  \right \rbrace = \nu\,\left \lbrace   y\, \in \,F\, : \quad |g(y)| > t \right \rbrace   $$

\end{definition}

\vskip 1cm

\begin{definition}\textbf{Symmetric rearrangement of a set of~$\R^d $}\\
\noindent Let~$\cal V$ be a measurable set of finite volume in~$\R^d$. Its symmetric rearrangement~${\cal V}^\star$ is the open centered ball, the volume of which agrees with~$\cal V$:

$$
{\cal V}^\star =\left \lbrace x \,\in\,\R^d \,   \bigg| \, {\cal V}ol_{{\cal B}_d}\cdot|x|^d < {\cal V}ol_{\cal V}  \right \rbrace
$$

\noindent where~${\cal V}_{{\cal B}_d}$ denotes the volume of the unit ball of~$\R^d$ :

$$
{\cal V}ol_{{\cal B}_d}= \displaystyle \frac {\pi^{\frac{d}{2}}}{\Gamma \left (\frac{d}{2}+1 \right)}
$$

\end{definition}
\vskip 1cm

\begin{definition}\textbf{Symmetric decreasing rearrangement}\\
Let us consider a real-valued, non-negative, measurable function $f$, defined on~$\R^d  $, which vanishes at infinity, in the sense where

$$
\text{mes}\, \left (  \left \lbrace x \,\in\,\R^d \,   \bigg| \,  f(x)  > t \right \rbrace  \right) < + \infty \quad \forall \,t > 0
$$

\noindent The symmetric decreasing rearrangement of~$f$ is the function $f^\star $, defined on~$\R^d $, positive, measurable, such that, for any~$x$ ~in $\R^d $:

$$ f^\star(x) = \displaystyle\int_0^{+\infty} \mathbbm{1}_{\{y: f(y)>t\}^\star}(x) \, dt$$

\noindent where~\mbox{$\mathbbm{1}_{\{y: f(y)>t\}^\star}$} denotes the characteristic function of the set~\mbox{$ {\{y: f(y)>t\}^\star}$}.

\end{definition}


\vskip 1cm

\begin{definition}\textbf{Decreasing rearrangement}\\

\noindent   If one denotes by~$f$ a real-valued, measurable function, $\mu-$finite a.e., the (scalar) decreasing rearrangement~$f^\#$, is the positive-valued function, defined, for any strictly positive number~$t$, through:

$$ f^\#(t)=   \displaystyle \inf_{\lambda > 0}\,  \left \lbrace \lambda \, \big |\, \mu_f(\lambda) \leq t\right \rbrace$$

\noindent with the convention:

$$  \displaystyle \inf \emptyset = + \infty$$

\end{definition}

\vskip 1cm

\begin{example}\textbf{Distribution function and rearrangement of a simple function}\\

\noindent Let us denote by $n$ a natural integer, and $a_1$, .., $a_n$ real numbers such that:

$$a_1< \hdots <   a_n$$

\noindent  $I_1$, ..., $I_n$ are pairwise disjoint intervals of~$\R $. We define the function $f$, from $\R$ to $\R$, such that:

$$f= \displaystyle \sum_{i=1}^n a_i\,\1_{I_i}$$

\noindent where, for~$1 \leq i\leq n_f$,  $\mathbbm{1}_{I_{i }}$ is the characteristic function of~$I_{i }$:

$$\forall\,x \,\in\,\R^d \, : \, \mathbbm{1}_{I_{i }}(x)= \left \lbrace \begin{array}{ccc} 1 & \text{if} & x\,\in\,I_{i } \\ 0 & \text{otherwise} \end{array} \right.$$

\noindent Then:\\

\begin{enumerate}

\item[$\rightsquigarrow$] For any $\lambda \geq a_n$:

$$\mu_f(\lambda)=0$$

\item[$\rightsquigarrow$] For any~$\lambda \,\in\,\left[ a_{n-1},a_n\right[ $:

$$ |f(x)| > \lambda   \Longleftrightarrow x\,\in\, I_n$$

\noindent which leads to:

$$\mu_f(\lambda)=\text{mes}\,I_n$$

\item[$\rightsquigarrow$] For any~$\lambda \,\in\,\left[ a_{n-2},a_{n-1}\right[ $:

$$  |f(x)| > \lambda   \Longleftrightarrow x\,\in\, I_{n-1} \cup I_n$$

\noindent One thus has:

$$\mu_f(\lambda)=\text{mes}\,I_{n-1} +\text{mes}\,I_n $$

\item[$\rightsquigarrow$] For any~$\lambda \,\in\,\left[ a_{i},a_{i+1}\right[ $:

$$x \,\in\,  \left \lbrace x \, \big |\, |f(x)| > \lambda \right \rbrace \Longleftrightarrow x\,\in\, I_{i} \cup \hdots \cup I_n$$

\noindent This yields:

$$\mu_f(\lambda)=\text{mes}\,I_{i+1} +\hdots+\text{mes}\,I_n $$

\item[$\rightsquigarrow$] If one sets $a_{n+1}=0$, one gets, by induction:

$$\mu_f =\displaystyle \sum_{i=1}^n \left \lbrace \displaystyle \sum_{k=i}^n \text{mes}\,I_{k}\right \rbrace \, \1_{[0,a_{i }[} $$

\noindent which can be written as:

$$\mu_f =\displaystyle \sum_{i=1}^n m_i \, \1_{[a_{i-1},a_{i }[} $$

\noindent with the convention:

$$a_0=0$$

\noindent and where, for any~$i$~in~$\left \lbrace 1,\hdots,n \right \rbrace $:

$$m_i=  \displaystyle \sum_{k=i}^n \text{mes}\,I_{k} $$

\noindent One also has:

$$ f^\# =\displaystyle \sum_{i=1}^n  a_i\, \1_{[m_{i-1},m_{i }[} $$

\noindent with the convention:

$$m_0=0$$

\end{enumerate}

\vskip 1cm

\noindent Figures 1, 2, 3 display the graphs of~$f$, $\mu_f$ et $f^\#$ in the case~$n=3$:

\begin{figure}[h!]
 \center{\psfig{height=7cm,width=9cm,angle=0,file=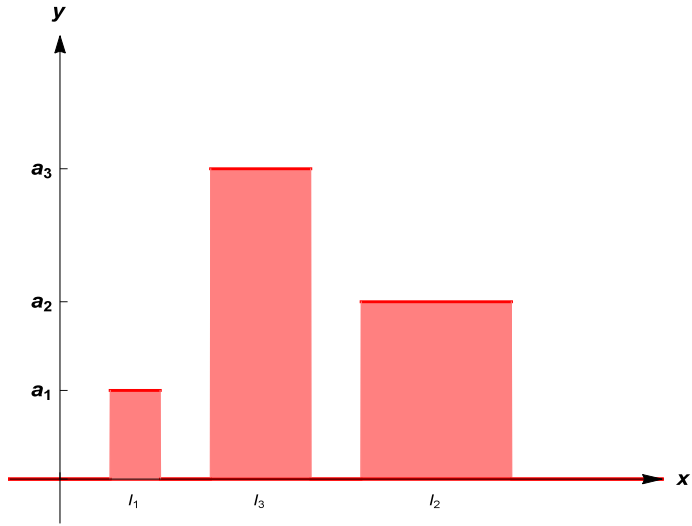}}
 \caption{The graph of the function~$f$ for~$n=3$.}

\end{figure}

\begin{figure}[h!]
 \center{\psfig{height=9cm,width=12cm,angle=0,file=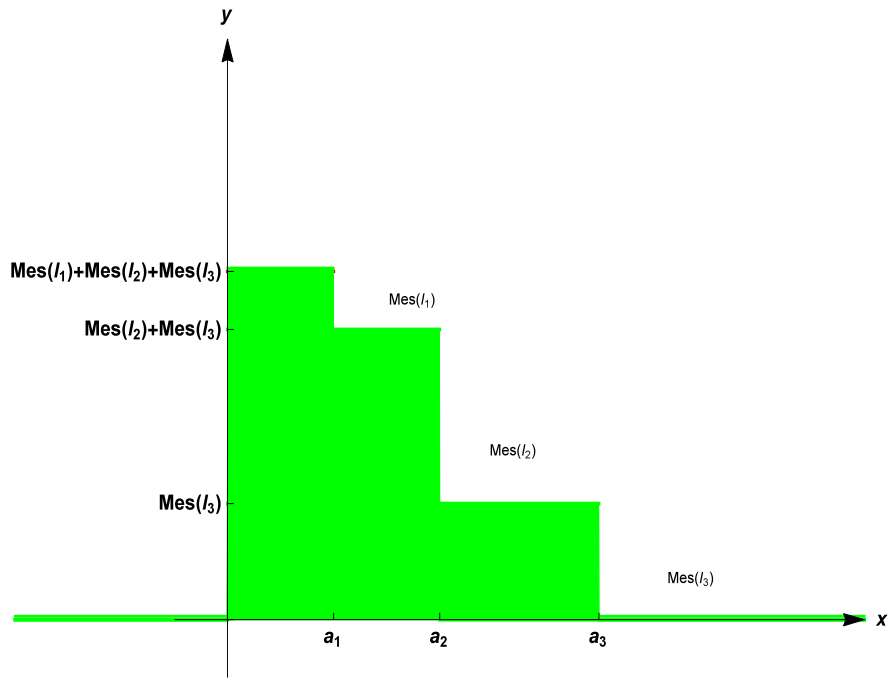}}
 \caption{The graph of the distribution function of~$f$ for~$n=3$.}

\end{figure}

\begin{figure}[h!]
 \center{\psfig{height=9cm,width=12cm,angle=0,file=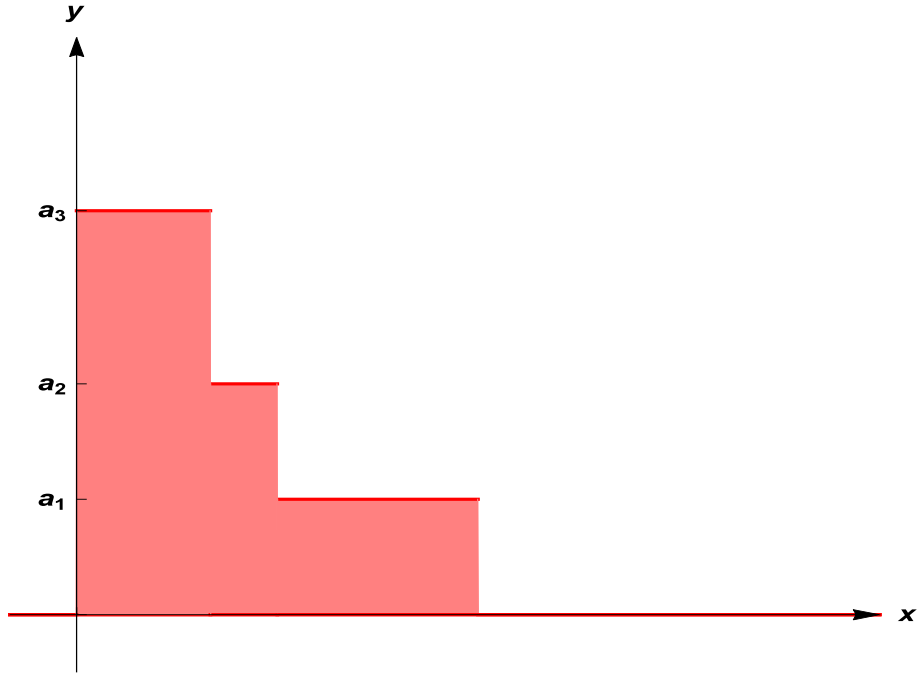}}
 \caption{The graph of the decreasing rearrangement~$f^\#$ for~$n=3$.}
\end{figure}

\end{example}

\vskip 1cm

\begin{pte}\textbf{Some properties of the decreasing rearrangement}\\

\noindent  Let us denote by~$f$ a real-valued, measurable function, $\mu-$finite a.e. Then:


\begin{enumerate}

\item[\emph{i}.] For any strictly positive number~$t$, and any positive number~$\lambda$:

$$f^\#(t) > \lambda \Longleftrightarrow  \mu_f(\lambda) >t$$

\item[\emph{ii}.] The function~$f^\#$  is non increasing.

\item[\emph{iii}.] If the distribution function~$\mu_f$ is strictly decreasing and continuous from~$I_{\mu_f} \subset \R$ into $\R^+$, then~$f^\#$ is the inverse function of~$\mu_f$ on~$\R^+$. Moreover, the decreasing rearrangement~$f^\#$ is right-continuous.\\


\end{enumerate}

\end{pte}

\vskip 1cm

\begin{proof}

\begin{enumerate}

\item[\emph{i}.] Let us denote by~$\lambda$ a positive number. Then, if~$\mu_f(\lambda) >t$, one gets, due to the fact that $\mu_f $ is non-increasing:

$$\lambda > \displaystyle \inf\, \left \lbrace \nu \,\big |\, \mu_f(\nu) \leq t  \right \rbrace $$

\noindent i.e.:

$$\lambda > f^\#(t)$$

\noindent Conversely, let us assume that~$\lambda > f^\#(t)$, i.e.:

$$\lambda > \displaystyle \inf\, \left \lbrace \nu \,\big |\, \mu_f(\nu) \leq t  \right \rbrace $$

\noindent The distribution function~$\mu_f$ being non-increasing, one deduces: $\mu_f(\lambda) >t$.

For any strictly positive number~$t$:

$$f^\#(t) > \lambda \Longleftrightarrow \mu_f(\lambda) >t$$

\item[\emph{ii}.] The function~$f^\#$  is decreasing; for~$t_1 \leq t_2$:

$$ \left \lbrace x \, \big |\,  |f(x) | > t_1\right \rbrace  \subset  \left \lbrace x \, \big |\,  |f(x) | > t_2\right \rbrace  $$

\noindent Thus:

$$ \text{mes}\left ( \left \lbrace x \, \big |\,  |f(x) | > t_1\right \rbrace \right) \leq  \text{mes} \left ( \left \lbrace x \, \big |\,  |f(x) | > t_2\right \rbrace \right)$$

\noindent i.e.:

$$ f^\#(t_1) \leq f^\#(t_2) $$

\item[\emph{iii}.]  The fact that the decreasing rearrangement~$f^\#$ is right-continuous follows from the fact that it is the distribution function of~$\mu_f$, with respect to the Lebesgue measure on~$\R^+$:

$$\forall  \,t \geq 0 \, : \quad
 f^\#(t)=\inf \, \left \lbrace \lambda \,\big |\, \mu_f (\lambda )\leq  t \right \rbrace =\sup \, \left \lbrace \lambda \,\big |\, \mu_f (\lambda )> t \right \rbrace =m_{\mu_f}(t)$$


\noindent Thus, its is right-continuous.     \\

\end{enumerate}

\end{proof}

\vskip 1cm

\begin{pte}\textbf{Hardy-Littlewood inequality}(One may refer to~\cite{BennettSharpley})\\

\noindent  Let us denote by~$f$ and~$g$ two measurable, real-valued functions defined on $\R$, which vanish at infinity. Then

$$\displaystyle    \int_{\R^d} f(x)\,g(x) \, dx \leq  \displaystyle    \int_{\R^+} f^\#(s)\,g^\#(s)  \,ds$$

\end{pte}

\vskip 1cm

\begin{proof}
The proof is made in the case of non-negative functions $f$ and $g$. Due to the monotone convergence theorem, one just needs to consider the case of a simple function~$f$, of the form:

$$f=  \displaystyle\sum_{i=1}^{n_f} f_{i }\,\mathbbm{1}_{A_{i }}  $$

\noindent where~$n_f$ is a positive integer, and, for~$1 \leq i\leq n_f$, $f_{i }$ is a positive number; $A_{1 }$, .., $A_{n_f }$ are measurable sets such that:

$$A_1 \subset  A_2 \subset \hdots \subset  A_{n_f }$$

\noindent One has then:

$$f^\#=\displaystyle\sum_{i=1}^{n_f} f_{i }\,\mathbbm{1}_{[0, \mu (A_{i })]}  $$

\noindent which leads to:

$$\begin{array}{ccc}\displaystyle    \int_{\R^d} f(x)\,g(x) \, dx &=&
\displaystyle\sum_{i=1}^{n_f} f_{i }\,\displaystyle \int_{A_i} g(x) \, dx\\
& \leq &
\displaystyle\sum_{i=1}^{n_f} f_{i }\, \displaystyle\int_{0}^{\mu (A_{i })} g^\#(t) \, dt\\
& \leq &
 \displaystyle\int_{\R^+} \displaystyle\sum_{i=1}^{n_f} f_{i }\,\,\mathbbm{1}_{[0, \mu (A_{i })]}(t)\, g^\#(t) \, dt\\
&= &
 \displaystyle\int_{\R^+} f^\#(t)\, g^\#(t) \, dt\\
\end{array}$$

\end{proof}

\vskip 1cm

\begin{proposition}\textbf{Contractive properties of the non-increasing rearrangement}\\
\label{resumpropstar}
\noindent Let us consider two real-valued, measurable functions $f$ and $g$, defined on $\R^d$, which vanishes at infinity. For any strictly positive number~$s$:
$$
  \|f^\#-g^\#\|_{L^p} \leq \|f-g\|_{L^p}
$$

\end{proposition}

\vskip 1cm

\begin{remark}
We would like to point out, as it appears in~\cite{LiebLoss1997}, that this specific property is a generalization of a theorem of G.~Chiti and C.~Pucci~\cite{ChitiPucci1979} and~\cite{CrandallTartar1980}. It however claimed to be proved for the first time in the 1986 paper~\cite{CroweZweibelRosenbloom1986}, which does not mention at all the work of ~G.~Chiti and C.~Pucci.

\end{remark}

\vskip 1cm

\begin{pte}\textbf{Rearrangement and algebric properties}\\

\noindent Let us denote by~$f$ and~$g$ two real-valued, measurable functions, $\mu-$finite a.e. Then, for all~$(t_1,t_2)~$ belonging to~$\R_+^2$:

$$ \left (f\,g\right)^\#(t_1+t_2) \leq f^\#(t_1)\,g^\#(t_2)  \quad, \quad  \left (f+g\right)^\#(t_1+t_2) \leq f^\#(t_1)+g^\#(t_2)$$

\end{pte}

\vskip 1cm

\begin{proof}

For any~$(t_1,t_2)\,\in\,\,\R_+^2\,$, one can legitimately assume that the quantities $f^\#(t_1) \,g^\#(t_2)$ and $f^\#(t_1) +g^\#(t_2)$ are respectively finite, since there is, otherwise, nothing to prove.\\

\noindent One has, first:

$$\left \lbrace x \, \big |\, |f(x)\,g(x) |> \lambda_1\,\lambda_2 \right \rbrace  \subset
\left \lbrace x \, \big |\, |f(x)|  > \lambda_1   \right \rbrace \cup \left \lbrace x \, \big |\, |g(x)|  > \lambda_2  \right \rbrace $$

\noindent since:\\

\noindent $\leftrightsquigarrow$ if $|f(x)\,g(x) |> \lambda_1\,\lambda_2$ et $|f(x)|  < \lambda_1 $, then, necessarily $|g(x)|  > \lambda_2 $ ;

\noindent $\leftrightsquigarrow$ if $|f(x)\,g(x) |> \lambda_1\,\lambda_2$ et $| g(x)|  < \lambda_2 $, then, necessarily  $|f(x)|  > \lambda_1 $.\\

\noindent Hence:

$$\text{mes} \, \left \lbrace x \, \big |\, |f(x)\,g(x) |> \lambda_1\,\lambda_2  \right \rbrace  \leq
\text{mes} \,\left \lbrace x \, \big |\, |f(x)|  > \lambda_1   \right \rbrace \cup \left \lbrace x \, \big |\, |g(x)|  > \lambda_2   \right \rbrace$$

\noindent i.e.:

$$\text{mes} \, \left \lbrace x \, \big |\, |f(x)\,g(x) |> \lambda_1\,\lambda_2  \right \rbrace  \leq
\text{mes} \, \left \lbrace x \, \big |\, |f(x)|  > \lambda_1   \right \rbrace +\text{mes} \, \left \lbrace x \, \big |\, |g(x)|  > \lambda_2   \right \rbrace$$

\noindent One has then:

$$\mu_{f\,g}(\lambda_1\,\lambda_2) \leq \mu_f(\lambda_1)+\mu_f(\lambda_2)$$

\noindent which, again due to the fact that the rearrangement is non-increasing, leads to:

$$\left (f\,g\right)^\# \left(\mu_f(\lambda_1)+\mu_f(\lambda_2) \right)  \leq \left (f\,g\right)^\# \left( \mu_{f\,g}(\lambda_1\,\lambda_2) \right)$$

\noindent i.e.:

$$\left (f\,g\right)^\# \left(\mu_f(\lambda_1)+\mu_f(\lambda_2) \right)  \leq  \lambda_1\,\lambda_2 $$

\noindent If one sets $t_1=\mu_f(\lambda_1)$, $t_2=\mu_g(\lambda_2)$, or, equivalently:

$$ \lambda_1= f^\#(t_1) \quad, \quad \lambda_2= g^\#(t_2) $$

\noindent one gets:

$$\left (f\,g\right)^\# \left(t_1+t_2 \right)  \leq  f^\#(t_1) \,g^\#(t_2) $$

\noindent Let us notice then that, for any couple of positive numbers~$(u,v)$, and for any~\mbox{$(a,b)\,\in\,\R^2$}, one has:

$$|a+b|> u+v \Rightarrow |a |> u \quad \text{or} \quad | b|>  v$$

\noindent since:

$$\left (|a | \leq u \quad \text{and} \quad |b|\leq v  \right)\Rightarrow |a+b | \leq |a |+|b | \leq u+v$$


$$  f^\#(t_1)+g^\#(t_2) < |f(x)+g(x) | <|f(x)|+|g(x) | $$
\noindent One has then the natural embedding


$$\left \lbrace x \, \big |\, |f(x)+g(x) |> f^\#(t_1)+g^\#(t_2)  \right \rbrace
\subset \left \lbrace x \, \big |\, |f(x)  |> f^\#(t_1)  \right \rbrace \cup \left \lbrace x \, \big |\, |g(x)  |> g^\#(t_2)  \right \rbrace$$



\noindent Thus:

$$\begin{array}{ccc}
\mu_{f+g}\left ( f^\#(t_1)+g^\#(t_2)  \right)&=&\text{mes}\, \left \lbrace x \, \big |\, |f(x)+g(x) |> f^\#(t_1)+g^\#(t_2)  \right \rbrace \\
& \leq &\text{mes}\, \left \lbrace x \, \big |\, |f(x)  |> f^\#(t_1)  \right \rbrace+\text{mes}\, \left \lbrace x \, \big |\, | g(x) |>  g^\#(t_2)  \right \rbrace  \\
& =&\mu_f  \left ( f^\#(t_1)  \right )+\mu_g  \left (  g^\#(t_2)  \right ) \\
& \leq & t_1+t_2  \\
\end{array}
$$

\noindent Since the rearrangement $\left (f+g\right)^\# $ is decreasing, it yields:

$$\left (f+g\right)^\#(t_1+t_2) \leq \left (f+g\right)^\#\left ( \mu_{f+g}\left ( f^\#(t_1)+g^\#(t_2)  \right) \right)
=(f+g )^\#\left ( \mu_{f+g}\left ( f^\#(t_1)+g^\#(t_2)  \right)\right)
=  f^\#(t_1)+g^\#(t_2) $$

\end{proof}

\vskip 1cm
\begin{definition}\textbf{Maximal function}\\

\noindent Let us denote by~$f$ a function defined on~$\R$, real-valued, measurable, finite a.e. We introduce the maximal function $f^{\star \star}$, defined, for any strictly positive number~$t$, by:

$$f^{\star \star} (t) = \displaystyle \frac{1}{t}\,\displaystyle \int_0^{t}  f^{\#} (s)\,ds$$

\end{definition}

\vskip 1cm

\begin{pte}
\label{PteMaxFunction}

\noindent Let us denote by~$f$ a real-valued, measurable function, $\mu-$finite a.e. The maximal function~$f^{\star \star}$ is non-increasing on~$\R_+^\star$. Moreover, for any strictly positive number~$t$:\\






$$f^{\#} (t) \leq f^{\star \star} (t)$$

\end{pte}

\vskip 1cm

\begin{proof}

\noindent The maximal function~$f^{\star \star}$ is non-increasing on~$\R_+^\star$ because, for any set of strictly positive numbers~$(t_1,t_2)$ such that $t_1 \leq t_2$:

$$\begin{array}{ccc}{f^{\star \star}}  (t_2) &=&\displaystyle \frac{1}{t_2}\,\displaystyle \int_0^{t_2} f^\#(s)  \,ds\\
 &=&\displaystyle \frac{1}{t_2}\,\displaystyle \int_0^{t_1} f^\#(s)  \,ds+\displaystyle \frac{1}{t_2}\,\displaystyle \int_{t_1}^{t_2} f^\#(s)  \,ds\\
 & \leq &\displaystyle \frac{1}{t_2}\,\displaystyle \int_0^{t_1} f^\#(s)  \,ds+\displaystyle \frac{1}{t_2}\,\displaystyle \int_{t_1}^{t_2} f^\#(t_1)  \,ds\\
 & =&\displaystyle \frac{1}{t_2}\,\displaystyle \int_0^{t_1} f^\#(s)  \,ds+\displaystyle \frac{(t_2-t_1)\,f^\#(t_1)}{t_2} \\
 & =&\displaystyle \frac{t_1}{t_2}\,{f^{\star \star}}  (t_1)+ \displaystyle \frac{ t_2-t_1  }{t_1\,t_2}\,\displaystyle \int_0^{t_1} f^\#(t_1)  \,ds \\
 & \leq&\displaystyle \frac{t_1}{t_2}\,{f^{\star \star}}  (t_1)+ \displaystyle \frac{ t_2-t_1  }{t_1\,t_2}\,\displaystyle \int_0^{t_1} f^\#(s)  \,ds \\
 & =&\displaystyle \frac{t_1}{t_2}\,{f^{\star \star}}  (t_1)+ \displaystyle \frac{ t_2-t_1  }{t_2}\,  f^{\star \star}   (t_1) \,ds \\
 & =& f^{\star \star}   (t_1)  \\
\end{array} $$

\noindent One can also note that, for any strictly positive number~$t$:

 $$\begin{array}{ccc}
 \displaystyle\frac{d   }{dt}{f^{\star \star}}(t) &=&\displaystyle \frac{f^\#(t)}{t}- \displaystyle \frac{1}{t^2}\,\displaystyle \int_0^{t} f^\#(s)  \,ds\\
 &=&\displaystyle \frac{1}{t^2}\,\left \lbrace t\,f^\#(t)-\displaystyle \int_0^{t} f^\#(s)  \,ds \right \rbrace \\
 &\leq &\displaystyle \frac{1}{t^2}\,\left \lbrace t\,f^\#(t)-f^\#(t) \,\displaystyle \int_0^{t}ds \right \rbrace \\
& =0
 \end{array} $$

\noindent since the rearrangement~$f^\#$ is non-increasing, which yields:

  $$ \displaystyle\frac{d   }{dt}{f^{\star \star}}  (t) \leq 0$$







\noindent    Finally, one has obviously:

$$\begin{array}{ccc}f^{\star \star} (t) &= &
  \displaystyle \frac{1}{t}\,\displaystyle \int_0^{t} f^\#  (s)\,ds\\
&\geq &
  \displaystyle \frac{1}{t}\,\displaystyle \int_0^{t} f^\#  (t)\,ds\\
&  = &
  \displaystyle \frac{1}{t}\,t\, f^\#  (t) \\
&  = & f^\#  (t) \\
\end{array}$$

\end{proof}

\vskip 1cm

\begin{pte}
\noindent The maximal function is sub-additive: if~$f$ and~$g$ are real-valued, measurable function, $\mu-$finite a.e., then:

$$\forall \,t>0\,: \quad   (f+g)^{\star \star} (t) \leq f^{\star \star} (t) +g^{\star \star} (t) $$

\end{pte}

\vskip 1cm

\begin{proof}


\noindent Since:

$$\forall\,(t_1,t_2)\,\in\,\,\R_+^2\, : \quad  \left (f+g\right)^\#(t_1+t_2) \leq f^\#(t_1)+g^\#(t_2)$$

\noindent one deduces, for any strictly positive number~$t$:

$$\begin{array}{ccc}(f+g)^{\star \star} (t) &= &
  \displaystyle \frac{1}{t}\,\displaystyle \int_0^{t} (f+g)^\#  (s)\,ds\\
  & \leq  &
  \displaystyle \frac{1}{t}\,\displaystyle \int_0^{t} \left \lbrace f^\#\left (\displaystyle \frac{s}{2}\right )+g^\#\left (\displaystyle \frac{s}{2}\right ) \right \rbrace \,ds\\
   & \leq  &
  \displaystyle \frac{1}{t}\,\displaystyle \int_0^{\frac{t}{2}} \left \lbrace f^\#(s)+g^\#(s ) \right \rbrace \,2\,ds\\
  & =  &
   \displaystyle \frac{2}{t}\,\displaystyle \int_0^{\frac{t}{2}}   f^\#(s)\,ds+
   \displaystyle \frac{2}{t}\,\displaystyle \int_0^{\frac{t}{2}}   g^\#(s)\,ds\\
       &  \leq  &
   \displaystyle \frac{2}{t}\,\displaystyle \int_0^{\frac{t}{2}}   f^\#\left (\displaystyle \frac{s}{2}\right )\,ds+
   \displaystyle \frac{2}{t}\,\displaystyle \int_0^{\frac{t}{2}}   g^\#\left (\displaystyle \frac{s}{2}\right )\,ds\\
    & =   &
   \displaystyle \frac{1}{t}\,\displaystyle \int_0^{t}   f^\#(s)\,ds+
   \displaystyle \frac{1}{t}\,\displaystyle \int_0^{t}   g^\#(s)\,ds\\
       & =   &
   f^{\star \star} (t)+
   g^{\star \star} (t) \\
\end{array}$$

 \noindent since the rearrangements~$f^\#$ and~$g^\#$ are both decreasing functions.

\end{proof}
\vskip 1cm

 \vskip 1cm

\begin{proposition}\textbf{From the decreasing rearrangement towards the symmetric one}\\

\noindent Let us denote by~$f$ a real-valued, measurable function, defined on~$\R^d$, which vanishes at infinity.\\

\noindent If~$\omega_d$ is the volume of the unit ball of~$\R^d$, then, for any~$x$ of $\R^d$:

$$f^\star(x)= f^\#(\omega_d\,|x|^d) $$

\noindent Especially:

$$\forall \,t \,\in\,\R\, :\quad f^\star(t)=\displaystyle \frac{1}{2}\, f^\#(|t|)$$




\end{proposition}

\vskip 1cm

\begin{pte}

\noindent Let us denote by~$f$ a real-valued, measurable function, defined on~$\R^d$, which vanishes at infinity. Then, for any~$t>0$ :

$$ \text{mes}\left ( \left \lbrace x\,\in\,\R^d \, \big |\, |f(x)| > t \right \rbrace \right) =
\text{mes}\left ( \left \lbrace x\,\in\,\R^d \, \big |\, |f^\star(x)| > t \right \rbrace \right) =
\text{mes}\left ( \left \lbrace s\,\in\,\R  \, \big |\, |f^\#( s)| > t \right \rbrace \right) $$

\end{pte}

\vskip 1cm

\begin{pte}\textbf{Rearrangement and dilatation}\\

\noindent Let us denote by~$f$ a real-valued, measurable function, defined on~$\R^d$, which vanishes at infinity. $\Lambda $~is a strictly positive real number. One has then:

$$\left (f (\Lambda\cdot) ) \right)^\#= f ^\# ( \Lambda^{ d}\,\cdot)$$

\noindent The decreasing rearrangement~$f^{\star  }$ is such that:

$$\left (f (\Lambda\,\cdot) ) \right)^\star= f ^\star ( \Lambda \,\cdot)$$

\noindent The maximal function~$f^{\star \star}$ ~satisfies:

$$\forall\,t> 0 \, : \quad \left (f (\Lambda\cdot)   \right)^{\star \star}(t)=f^{\star \star}(\Lambda^{ d}\,t)$$

\end{pte}

\vskip 1cm

\begin{proof}

$$ \begin{array}{ccc}\left (f  (\Lambda\,\cdot)\right)^\# (t)
&=& \displaystyle \inf_{\sigma > 0}\,  \left \lbrace \sigma \, \big |\, m(\sigma, f(\Lambda\cdot)) \leq  t\right \rbrace\\
\end{array}$$

\noindent For any strictly positive number $\sigma $:

$$\begin{array}{ccc}  m(\sigma , f (\Lambda\cdot) )&= & \text{mes}\left ( \left \lbrace x \, \big |\, |f (\Lambda\,x)|  > \sigma\right \rbrace \right)\\
&=&\displaystyle \int_{\R^d}  \1_{|f (\Lambda\,x)|>\sigma}\, dx \\
&=&\displaystyle \int_{\R^d}  \1_{|f ( x)|>\sigma} \,\Lambda^{-d}\,dx \\
&=&\Lambda^{-d}\,m(\sigma , f   )\\
\end{array}$$

\noindent Thus:

$$ \begin{array}{ccc} \left (f  (\Lambda\,\cdot)\right)^\# (t)
&=& \displaystyle \inf_{\sigma > 0}\,  \left \lbrace \sigma \, \big |\, m(\sigma, f(\Lambda\cdot)) \leq  t\right \rbrace\\
&=& \displaystyle \inf_{\sigma > 0}\,  \left \lbrace \sigma \, \big |\, \Lambda^{-d}\,m(\sigma , f   ) \leq  t\right \rbrace\\
&=& \displaystyle \inf_{\sigma > 0}\,  \left \lbrace \sigma \, \big |\,  m(\sigma , f   ) \leq  \Lambda^{ d}\,t\right \rbrace\\
&=& f ^\# ( \Lambda^{ d}\,t)
\end{array}$$

\noindent Hence:

$$f^\star(\Lambda\, x)= f^\#(\omega_d\,\Lambda^d\,|x|^d)
 $$

\noindent Or:

$$  \left (f  (\Lambda\,\cdot)\right)^\star (x)=\left (f  (\Lambda\,\cdot)\right)^\#(\omega_d\,|x|^d)
= f^\#(\omega_d\,\Lambda^d\,|x|^d) = f^\star(\Lambda\, x) $$

\noindent Also, for any strictly positive number~$t$:

$$ \displaystyle \int_0^{t} \left (f (\Lambda\cdot)   \right)^\# (s)\,ds = \displaystyle \int_0^{t} f ^\# ( \Lambda^{ d}\,s)\,ds
= \displaystyle \int_0^{\Lambda^{ d}\,t} f ^\# (  s)\,\Lambda^{- d}\,ds$$

$$ \displaystyle \frac{1}{t}\,\displaystyle \int_0^{t} \left (f (\Lambda\cdot)   \right)^\# (s)\,ds =
\displaystyle \frac{1}{\Lambda^{  d}\,t}\, \displaystyle \int_0^{\Lambda^{ d}\,t} f ^\# (  s)\, ds
=f^{\star \star}(\Lambda^{ d}\,t)$$

\noindent i.e.:

$$\left (f (\Lambda\cdot) ) \right)^{\star \star}(t) =f^{\star \star}(\Lambda^{ d}\,t)
$$

\noindent Moreover:

$$\begin{array}{ccc}
 m(\sigma, \Lambda\,f)&=&\text{mes}\left ( \left \lbrace x \, \big |\, |f (x)|  > u^{\frac{1}{p}}\right \rbrace \right)\\
\end{array}
$$

$$ f^\#(t)=   \displaystyle \inf_{\sigma > 0}\,  \left \lbrace \sigma \, \big |\, m(\sigma, f) \leq t\right \rbrace$$

$$ \begin{array}{ccc}(f\,g)^\#( t)= m( t,f\,g)=  \displaystyle \inf_{\sigma > 0}\,  \left \lbrace \sigma \, \big |\, m(\sigma,  f\,g) \leq  t\right \rbrace\\
\end{array}$$

\noindent Also:

$$\begin{array}{ccc}  m(u^{\frac{1}{p}}, f  )&= & \text{mes}\left ( \left \lbrace x \, \big |\, |f (x)|  > u^{\frac{1}{p}}\right \rbrace \right)\\
&=&\text{mes}\left ( \left \lbrace x \, \big |\, |f (x)|^p  > u\right \rbrace \right) \\
&=&\displaystyle \inf_{s >0}  \,\left \lbrace s \, \big |\,\text{mes}\left ( \left \lbrace x \, \big |\, |f(x)|^p > s\right \rbrace \right)\leq u \right \rbrace\\
&=&\displaystyle \inf_{s >0}  \,\left \lbrace s \, \big |\,\text{mes}\left ( \left \lbrace x \, \big |\, |f(x)|  > s^{\frac{1}{p}}\right \rbrace \right)\leq u \right \rbrace\\
&=&\displaystyle \inf_{s >0}  \,\left \lbrace s^p \, \big |\,
 \text{mes}\left ( \left \lbrace x \, \big |\, |f(x)|  > s  \right \rbrace \right)\leq u  \right \rbrace\\
&= &m^p(u , f  )\\
&=&\left (f^\#\right)^p(u)
\end{array}$$

\noindent For any strictly positive number~$\Lambda$:

$$\begin{array}{ccc}  m(t , f (\Lambda\cdot) )&= & \text{mes}\left ( \left \lbrace x \, \big |\, |f (\Lambda\,x)|  > t\right \rbrace \right)\\
&=&\displaystyle \int_{\R^d}  \1_{|f (\Lambda\,x)|>t}\, dx \\
&=&\displaystyle \int_{\R^d}  \1_{|f ( x)|>t} \,\Lambda^{-d}\,dx \\
&=&\Lambda^{-d}\,m(t , f   )\\
\end{array}$$

\end{proof}

\vskip 1cm

\begin{pte}\textbf{Rearrangement and invariance of the~$L^p$ norm}\\

\noindent Let us consider $p\geq 2$; we denote by~$f$ a real-valued, measurable function, defined on~$\R^d$, which vanishes at infinity.
\noindent If the function~$f$ belongs to~$L^p \left ( \R^d \right)$, the symmetric decreasing rearrangement~$f^\star $ belongs to~$L^p \left ( \R^d \right)$, the decreasing rearrangement~$f^\# $ belongs to~$L^p \left ( \R \right)$, and, for any~$t > 0$, one has:

$$\left \|f^\star  \right \|_{L^p \left ( \R^d \right)} =\|f^\#    \|_{L^p \left ( \R  \right)}=\left \|f  \right \|_{L^p \left ( \R^d \right)}
$$

 $$
\left \|f  \right \|_{L^p \left ( \R^d \right)}^p= \displaystyle \int_0^{+\infty}   \left( u^\frac{1}{p}\, f^\# (u)\right)^p   \,\displaystyle \frac{du}{u}$$

\end{pte}

\vskip 1cm

\begin{proof}

\noindent The Fubini theorem leads to:

$$\begin{array}{ccc}
\left \| f \right \|^p_{L^p \left ( \R^d \right)}
 &=& \displaystyle\int_{\R^d} \left (f(x) \right )^p\,dx \\
  &=& \displaystyle\int_{\R^d} \left \lbrace  \displaystyle \int_0^{+\infty} \mathbbm{1}_{ (f(x)  )^p > t} \,dt \right \rbrace \,dx \\
  &=& \displaystyle \int_0^{+\infty}   \left \lbrace  \displaystyle \int_{\R^d}  \mathbbm{1}_{ (f(x)  )^p > t} \,dx \right \rbrace\,dt \\
  &=& \displaystyle \int_0^{+\infty}   \text{mes}\,\left \lbrace x\, \bigg| \,  (f(x)  )^p > t    \right \rbrace \,dt \\
  &=& \displaystyle \int_0^{+\infty}   \text{mes}\,\left \lbrace x\, \bigg| \,   f(x)    > t    \right \rbrace \,p\,t^{p-1}\,dt \\
   &=& \displaystyle \int_0^{+\infty}   \mu_f(t)\,p\,t^{p-1}\,dt \\
&=&\left \|f^\star  \right \|^p_{L^p \left ( \R^d \right)}
\end{array}$$

\noindent Also:

 $$\begin{array}{ccc}
\left \| f \right \|^p_{L^p \left ( \R^d \right)}
 &=& \displaystyle\int_{\R^d} |f(x)|^p\,dx \\
  &=& \displaystyle\int_{\R^d} \left \lbrace \displaystyle \int_0^{|f (x)|} p\,t^{p-1 }\, dt  \right \rbrace \,dx\\
  &=& \displaystyle \int_0^{+\infty} p\,t^{p-1 }\,\left \lbrace \displaystyle\int_{\R^d} \1_{|f (x)|>t} dx\,  \right \rbrace \,dt\\
  &=& \displaystyle \int_0^{+\infty} p\,t^{p-1 }\,\text{mes}\left ( \left \lbrace x \, \big |\, |f(x)| > t\right \rbrace \right)   \,dt\\
  &=& \displaystyle \int_0^{+\infty} p\,t^{p-1 }\,
\displaystyle \inf_{\sigma  > 0}\,  \left \lbrace \sigma  \, \big |\,  m\left (\sigma, f^\star \right) \leq t \right \rbrace   \,dt\\
 &=& \displaystyle \int_0^{+\infty} p\,t^{p-1 }\,
\displaystyle \inf_{\sigma  > 0}\,  \left \lbrace \sigma^p  \, \big |\,  m \left (\sigma^{\frac{1}{p}}, f^\star \right) \leq t  \right \rbrace   \,dt\\
  &=& \displaystyle \int_0^{+\infty} p\,t^{p-1 }\,m(t,f)   \,dt\\
&=& \displaystyle \int_0^{+\infty} p\,t^{p-1 }\,m(t,f^\star)   \,dt\\
&=& \displaystyle \int_0^{+\infty} p\,t^{p-1 }\,m\left ((t^p)^\frac{1}{p}  , f^\star \right)   \,dt\\
&=& \displaystyle \int_0^{+\infty}   m\left (u^\frac{1}{p}, f^\star \right)   \,du\\
&=& \displaystyle \int_0^{+\infty}   \left (f^\#\right)^p(u)   \,du\\
&=&\|f^\#    \|^p_{L^p \left ( \R  \right)}\\
&=& \displaystyle \int_0^{+\infty}   u\,\left (f^\#(u) \right)^p   \,\displaystyle \frac{du}{u}\\
&=& \displaystyle \int_0^{+\infty}   \left( u^\frac{1}{p}\, f^\# (u)\right)^p   \,\displaystyle \frac{du}{u}\\
\end{array}$$

\noindent An alternate proof of the relations between those norms can be made using the monotone convergence theorem. This way, one just needs to consider the case of a simple function~$f$, of the form:

$$f=  \displaystyle\sum_{i=1}^{n_f} f_{i }\,\mathbbm{1}_{A_{i }}  $$

\noindent where~$n_f$ is a positive integer, and, for~$1 \leq i\leq n_f$, $f_{i }$ is a positive number ; $A_{1 }$, .., $A_{n_f }$ are measurable sets such that:

$$A_1 \subset  A_2 \subset \hdots \subset  A_{n_f }$$

\noindent One has then:

$$f^\#=\displaystyle\sum_{i=1}^{n_f} f_{i }\,\mathbbm{1}_{[0, \mu (A_{i })]}  $$

\noindent For~$1 \leq i\leq n_f$, the equimeasurability of $f_i$ and $f_i^\#$ can be written, for any positive number $\lambda$, as:

 $$\mu_{f_i}(\lambda)=\mu_{f_i^\#}(\lambda)$$

 \noindent which leads to:

 $$\mu_{f }(\lambda)=\mu_{f^\#}(\lambda)$$

 \noindent and:

 $$\mu_{f^p }(\lambda)=\mu_{f }(\lambda^{\frac{1}{p}}) =\mu_{f^\#}(\lambda^{\frac{1}{p}}) =\mu_{(f^\#)^p}(\lambda )$$

\noindent i.e.:

$$ \mu \left ( \left \lbrace x \, \big |\, |f^p(x)| > \lambda\right \rbrace \right)= \mu \left ( \left \lbrace x \, \big |\, |(f^\#)^p(x)| > \lambda\right \rbrace \right)$$

 \vskip 1cm

\end{proof}

\vskip 1cm

\noindent Starting from the above latter property, Lorentz spaces can be introduced very naturally:

\vskip 1cm

\begin{definition}
\noindent Let~$p$ and $q$ denote two strictly positive numbers such that~$p>1$,~$q \geq1$. The Lorentz space~$L^{p,q}(\R^d)$ is defined as the set of real-valued, measurable functions $f$, defined on~~$\R^d$, such that:

$$\left \|f\right\|_{L^{p,q}(\R^d)}= \left ( \displaystyle \int_0^{+\infty} \left (t^{\frac{1}{p}} \,f^\#\ (t)\right)^q \,\displaystyle \frac{dt}{t} \right)^{\frac{1}{q}}
< + \infty$$

\end{definition}

\vskip 1cm

\begin{remark}
It is interesting to note that one can also define:\\

\begin{enumerate}
\item[\emph{i}.] the Lorentz space~$L^{p,\infty}(\R^d)$ as the set of real-valued, measurable functions $f$, defined on~$\R^d$, such that:

$$\left \|f\right\|_{L^{p,\infty}(\R^d)}=  \displaystyle \sup_{t>0}\,   t^{\frac{1}{p}} \,f^\#\ (t)< + \infty$$

\item[\emph{ii}.] the Lorentz space~$L^{\infty,\infty}(\R^d)$ as the set of real-valued, measurable functions $f$, defined on~$\R^d$,such that:

$$\left \|f\right\|_{L^{\infty,\infty}(\R^d)}=  \displaystyle \sup_{t>0}\,    f^\#\ (t)< + \infty$$

\end{enumerate}

\end{remark}
\vskip 1cm

\begin{remark}
\noindent It is clear that:

$$L^{p, p}(\R^d)=L^{p }(\R^d)$$

\noindent since the $L^p$~-norm is kept invariant by the rearrangement.


\end{remark}

\vskip 1cm

\begin{remark}

\noindent  Lorentz spaces can be considered as " thinner" spaces than the Lebesgue ones; they make it possible to
detect logarithmic correction, which can not be done with the classical $L^p$ spaces.

\vskip 1cm
\begin{example}

\noindent Let us consider, on the unit ball of $\R$, functions of the form:

$$f_{\alpha,\beta} \, : \, t \neq 0 \mapsto \displaystyle \frac{\left (\ln \frac{1}{|t|}\right)^\beta }{|t|^\alpha}
=\displaystyle \frac{  (-\ln  |t|  )^\beta}{|t|^\alpha}$$

$$
\begin{array}{ccc}\left \|  f_{\alpha,\beta }  \right\|_{{\cal B}_{0,1} (\R )}^p&=&
\displaystyle \int_{{\cal B}_{0,1} (\R )} f^p_{\alpha,\beta } (t)\,dt\\
&=&
\displaystyle \int_{{\cal B}_{0,1} (\R )} \displaystyle \frac{\left (\ln \frac{1}{|t|}\right)^{\beta \,p} }{|t|^{\alpha\,p}}\,dt\\
&=&
\displaystyle \int_{{\cal B}_{0,1} (\R )} \displaystyle \frac{\left (-\ln  |t| \right)^{\beta \,p} }{t^{\alpha\,p}}\,dt\\
\end{array}
$$

\noindent At stake are Bertrand integrals, of the form
$$\displaystyle \int_0^1  \displaystyle \frac{dt}{t^{\alpha_0}\,\left (-\ln t \right)^{\beta_0} }$$

\noindent When~$ \alpha_0<1$, no difference can be seen for distinct values of the parameter~$\beta$. It is not the case if one consider Lorentz norms, as it is illustrated in the following figures.



\begin{figure}[h!]
 \center{\psfig{height=8cm,width=8cm,angle=0,file=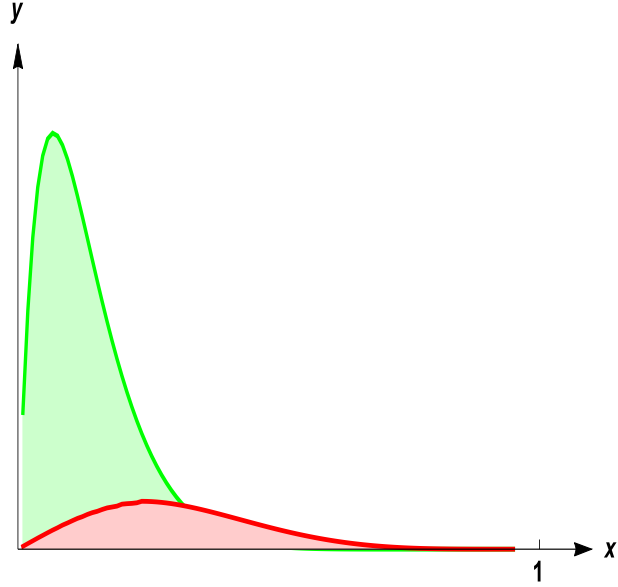} }
 \caption{The $\|\cdot\|_{L^{1,2}([0,1])}^2$ norms of the function $f_{-1,2}$ (in red) and $f_{-1,4}$ (in green).}

\end{figure}

\begin{figure}[h!]
 \center{\psfig{height=8cm,width=8cm,angle=0,file=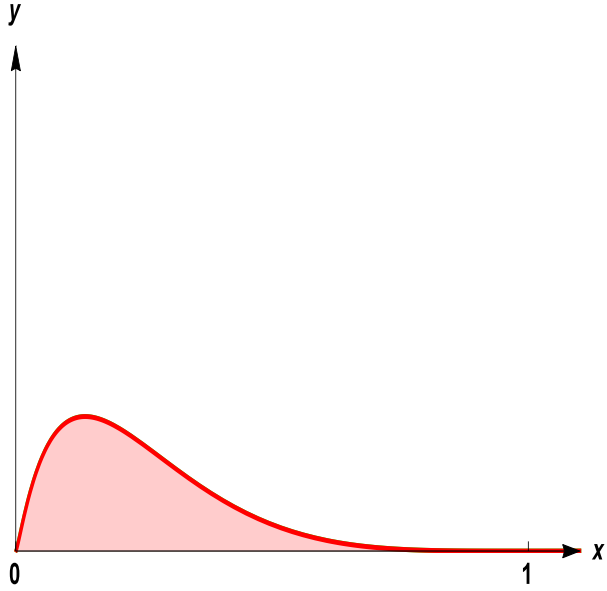} }
 \caption{The $\|\cdot\|_{L^{ 2}([0,1])}^2$ norms of the function $f_{-1,2}$ and $f_{-1,4}$, which are the same.}

\end{figure}

\end{example}

\vskip 1cm

 \end{remark}

\vskip 1cm


\begin{remark}

\noindent It is important to note that $\left \| \cdot \right\|_{L^{p,q}(\R^d)}$ is not a norm, since the triangle inequality does not hold, for, in most cases, one cannot have:

 $$\left (f+g \right)^\star \leq f^\star+g^\star$$

 \noindent $\left \| \cdot \right\|_{L^{p,q}(\R^d)}$ is just a quasi-norm.\\

 \noindent The space~$L^{p,q}(\R^d)$ is not, thus, a Banach space. In order to norm~$L^{p,q}(\R^d)$, one has to consider, thanks to the sub-additivity of the maximal function:

$$|||f  |||_{L^{p,q}(\R^d)}= \left ( \displaystyle \int_0^{+\infty} \left (t^{\frac{1}{p}} \,f^{\star \star} (t)\right)^q \,\displaystyle \frac{dt}{t} \right)^{\frac{1}{q}} \quad , \quad q < + \infty$$

\noindent The (obviously) vectorial space~$L^{p,q}(\R^d)$ is, thus, a complete metric space. It is obvious if one considers the mapping

 $$(f,g)\,\in\,L^{p,q}(\R^d) \times L^{p,q}(\R^d) \mapsto  ||| f-g |||_{L^{p,q}(\R^d)}$$

 \noindent which is a distance over~$L^{p,q}(\R^d)$. As for the mapping

 $$(f,g)\,\in\,L^{p,q}(\R^d) \times L^{p,q}(\R^d) \mapsto \left \| f-g \right\|_{L^{p,q}(\R^d)}$$

\noindent one can also prove that it is a distance over~$L^{p,q}(\R^d)$, which follows from the following comparison:

\end{remark}

\vskip 1cm

\begin{pte}\textbf{Comparison of the Lorentz norm and quasi-norm}\\
\noindent There exists a strictly positive constant~$C_{p,q}$ such that, for any~$f$ belonging to~$L^{p,q}(\R^d)$:

$$ |||  f  |||_{L^{p,q}(\R^d)} \leq   C_{p,q}\, ||f  ||_{L^{p,q}(\R^d)}  \leq   C_{p,q}\, |||f  |||_{L^{p,q}(\R^d)}$$

\end{pte}
\vskip 1cm

\begin{proof}

\noindent Let us recall:
$$|||f  |||_{L^{p,q}(\R^d)}= \left ( \displaystyle \int_0^{+\infty} \left (t^{\frac{1}{p}} \,f^{\star \star} (t)\right)^q \,\displaystyle \frac{dt}{t} \right)^{\frac{1}{q}} \quad,
\quad|| f  ||_{L^{p,q}(\R^d)}= \left ( \displaystyle \int_0^{+\infty} \left (t^{\frac{1}{p}} \,f^{\# } (t)\right)^q \,\displaystyle \frac{dt}{t} \right)^{\frac{1}{q}}  $$

\noindent Due to Hardy's generalized inequality, for any strictly positive number $T$, and any positive number~$ \alpha$ such that~\mbox{$\alpha < q-1$} :

$$\displaystyle \int_0^{T}  \left ( \displaystyle \frac{ \displaystyle \int_0^{t} f ^{\star}(u)\,du}{t}\right)^q \,t^\alpha\, dt \leq
\left (  \displaystyle \frac{q}{q-1-\alpha} \right)^q \,  \displaystyle \int_0^{T}  \left (f ^{\#}(t)\right)^q \,t^\alpha\,dt$$

\noindent Thus, for $\alpha= \displaystyle \frac{q}{p} $:

$$\displaystyle \int_0^{T}  \left (t^{\frac{1}{p}-1} \,  \displaystyle \int_0^{t} f ^{\star}(u)\,du \right)^q \,  dt \leq
  \left (  \displaystyle \frac{q}{q-1+\frac{q}{p}} \right)^q \, \displaystyle \int_0^{T} \left (t^{\frac{1}{p} } \,f ^{\#}(t)\right)^q  \,dt$$

\noindent i.e.:

$$\displaystyle \int_0^{T}  \left (t^{\frac{1}{p}} \,  f ^{\star\star}(t) \right)^q \,  dt \leq
  \left (  \displaystyle \frac{q}{q-1+\frac{q}{p}} \right)^q \, \displaystyle \int_0^{T} \left (t^{\frac{1}{p}} \,f ^{\#}(t)\right)^q  \,dt$$

\noindent which leads to:

$$\displaystyle \int_0^{+\infty}  \left (t^{\frac{1}{p}} \,  f ^{\star\star}(t) \right)^q \,  dt \leq
  \left (  \displaystyle \frac{q}{q-1+\frac{q}{p}} \right)^q \, \displaystyle \int_0^{+\infty} \left (t^{\frac{1}{p}} \,f ^{\#}(t)\right)^q  \,dt$$

\noindent and:

$$\left (\displaystyle \int_0^{+\infty}  \left (t^{\frac{1}{p}} \,  f ^{\star\star}(t) \right)^q \,  dt\right)^{\frac{1}{q}} \leq
   \displaystyle \frac{q}{q-1+\frac{q}{p}} \, \displaystyle \int_0^{+\infty} \left (t^{\frac{1}{p}} \,f ^{\#}(t)\right)^q  \,dt $$

\noindent i.e.:

$$|||f  |||_{L^{p,q}(\R^d)} \leq
   \displaystyle \frac{q}{q-1+\frac{q}{p}} \, || f   ||_{L^{p,q}(\R^d)} $$

\noindent or:

$$|||f  |||_{L^{p,q}(\R^d)} \leq
   \displaystyle \frac{p}{p+1-\frac{p}{q}} \, || f   ||_{L^{p,q}(\R^d)} =
 C_{p,q} \, || f   ||_{L^{p,q}(\R^d)} $$

\noindent The inequality

$$|| f   ||_{L^{p,q}(\R^d)} \leq||| f   |||_{L^{p,q}(\R^d)}   $$

\noindent follows from the comparison between~$f^\# $ and the maximal function~$f^{\star \star}$ seen in the above (see property (\ref{PteMaxFunction})).
\end{proof}

\vskip 1cm

\begin{pte}
\noindent The Lorentz space~$L^{p,q}(\R^d)$ is an homogeneous one ; for any strictly positive number~$\Lambda $, and any~$f$ belonging to $L^{p,q}(\R^d)$:

$$\left \|  f \left ( \Lambda \cdot \right )\right\|_{L^{p,q}(\R^d)}=\Lambda^{-\frac{d}{p} }\, \left \|f\right\|_{L^{p,q}(\R^d)}
\quad, \quad |||  f \left ( \Lambda \cdot \right )|||_{L^{p,q}(\R^d)}=\Lambda^{-d }\, |||f |||_{L^{p,q}(\R^d)}$$

\noindent Moreover:

$$\begin{array}{ccc}\left \|   \Lambda \,f  \right\|_{L^{p,q}(\R^d)}^q
&=&\Lambda \,\left \|  f  \right\|_{L^{p,q}(\R^d)}^q
\end{array}
\quad, \quad \begin{array}{ccc}|||    \Lambda \,f  ||| _{L^{p,q}(\R^d)}^q
&=&\Lambda \,|||   f  ||| _{L^{p,q}(\R^d)}^q
\end{array}
$$

\end{pte}

\vskip 1cm

\begin{proof}

$$\begin{array}{ccc}\left \|  f \left ( \Lambda \,\cdot \right )\right\|_{L^{p,q}(\R^d)}^q
&=&\displaystyle \int_0^{+\infty }  \left ( t^{\frac{1}{p}}\, \left (f \left ( \Lambda \cdot \right )\right )^\# (t) \right)^q\,\displaystyle \frac{dt}{t}\\
&=&\displaystyle \int_0^{+\infty }  \left ( t^{\frac{1}{p}}\,f^\#( \Lambda^{ d}\,t) \right)^q\,\displaystyle \frac{dt}{t}\\
&=&\displaystyle \int_0^{+\infty }  \left ( \Lambda^{ -\frac{d}{p}}\,t^{\frac{1}{p}}\,f^\#(  t) \right)^q\,\Lambda^{- d}\,\displaystyle \frac{dt}{\Lambda^{ -d}\,t}\\
&=&\Lambda^{ -\frac{d\,q} {p}}\,\displaystyle \int_0^{+\infty }  \left (  t^{\frac{1}{p}}\,f^\#(  t) \right)^p\, \displaystyle \frac{dt}{ t}\\
&=&\Lambda^{ -\frac{d\,q} {p}}\,\left \|  f  \right\|_{L^{p,q}(\R^d)}^q
\end{array}
$$

\noindent Hence:

$$\begin{array}{ccc}\left \|  f \left ( \Lambda \,\cdot \right )\right\|_{L^{p,q}(\R^d)}
&=&\Lambda^{- \frac{d}{p}}\,\left \|  f  \right\|_{L^{p,q}(\R^d)}
\end{array}
$$

\noindent For any strictly positive number~$ t$:

$$\begin{array}{ccc}\left (f (\Lambda\cdot) ) \right)^{\star \star} (t)
&= &\displaystyle \frac{1}{t}\,\displaystyle \int_0^{t}  \left (f (\Lambda\cdot)   \right)^\# (s)\,ds\\
&= &\displaystyle \frac{1}{t}\,\displaystyle \int_0^{t}  f^\#( \Lambda^{ d}\,s)\,ds\\
&= &\displaystyle \frac{1}{t}\,\displaystyle \int_0^{t} \Lambda^{-d}\,f^\# (s)\,ds\\
&=&\Lambda^{-d}\,f^{\star \star} (t) \\
\end{array}$$

\noindent Thus:
$$\left (f (\Lambda\cdot) ) \right)^{\star \star}=\Lambda^{-d}\,f^{\star \star}$$
\noindent and:

$$\begin{array}{ccc} |||  f \left ( \Lambda \,\cdot \right ) |||_{L^{p,q}(\R^d)}^q
&=&\displaystyle \int_0^{+\infty }  \left ( t^{\frac{1}{p}}\, \left (f \left ( \Lambda \cdot \right )\right )^{\star \star} (t) \right)^q\,\displaystyle \frac{dt}{t}\\
&=&\displaystyle \int_0^{+\infty }  \left ( t^{\frac{1}{p}}\,\Lambda^{ -d}\,f^{\star \star} (t)  \right)^q\,\displaystyle \frac{dt}{t}\\
&=&\displaystyle \int_0^{+\infty } \Lambda^{ - d\,q }\, \left (  t^{\frac{1}{p}}\,f^{\star \star} (  t) \right)^q\, \displaystyle \frac{dt}{  t}\\
&=&\Lambda^{ - d\,q }\,\displaystyle \int_0^{+\infty }  \left (  t^{\frac{1}{p}}\,f^{\star \star}  \right)^p\, \displaystyle \frac{dt}{ t}\\
&=&\Lambda^{-  d\,q }\,||| f |||_{L^{p,q}(\R^d)}^q
\end{array}
$$

\noindent Also:

$$\begin{array}{ccc}\left \|   \Lambda \,f  \right\|_{L^{p,q}(\R^d)}^q
&=&\Lambda \,\left \|  f  \right\|_{L^{p,q}(\R^d)}^q
\end{array}
$$

\end{proof}

\vskip 1cm

\begin{pte}

\noindent Let us consider $ 0 \leq p < + \infty$,  $ 0 \leq q_1 < + \infty$,  $ 0 \leq q_2 < + \infty$. Then:

$$ q_1 \leq q_2 \Rightarrow L^{p,q_1}(\R^d) \hookrightarrow L^{p,q_2}(\R^d)$$

\end{pte}
\vskip 1cm

\begin{proof}\cite{BennettSharpley}$\,$\\
\noindent Due to the fact that the rearrangement~$f^\#$ is non increasing, one gets, for any real positive number~$t$:

$$\begin{array}{ccc}
t^{\frac{1}{p}}\, f^\#(t)
&=& \left \lbrace \displaystyle \frac{p}{q_1}\,\displaystyle \int_0^t \left (s^{\frac{1}{p}}\, f^\#(t) \right)^{q_1} \,\displaystyle \frac{ds}{s}  \right \rbrace^{\frac{1}{q_1}} \\
&\leq & \left \lbrace \displaystyle \frac{p}{q_1}\,\displaystyle \int_0^t \left (s^{\frac{1}{p}}\, f^\#(s) \right)^{q_1} \,\displaystyle \frac{ds}{s}  \right \rbrace^{\frac{1}{q_1}} \\
&\leq & \left (\displaystyle \frac{p}{q_1}\right)^{\frac{1}{q_1}}\,  \| f \|{L^{p,q_1}(\R^d)} \\
\end{array}
$$

\noindent $\rightsquigarrow$ If~$q_1=+\infty$, one gets:

$$\begin{array}{ccc} \| f \|{L^{p,\infty}(\R^d)}
&\leq & \left (\displaystyle \frac{p}{q_1}\right)^{\frac{1}{q_1}}\,  \| f \|{L^{p,q_1}(\R^d)} \\
\end{array}
$$

\noindent $\rightsquigarrow$ If~$q_1<+\infty$, one gets:

$$\begin{array}{ccc} \| f \|{L^{p,q_2}(\R^d)}
&= & \left \lbrace \displaystyle \int_0^t \left (s^{\frac{1}{p}}\, f^\#(t) \right)^{q_2} \,\displaystyle \frac{ds}{s}  \right \rbrace^{\frac{1}{q_2}}\\ \\
&= & \left \lbrace \displaystyle \int_0^t \left (s^{\frac{1}{p}}\, f^\#(t) \right)^{q_2-q_1+q_1} \,\displaystyle \frac{ds}{s}  \right \rbrace^{\frac{1}{q_2}}\\ \\
&\leq & \| f \|^{\frac{q_2-q_1}{q_2}}_{L^{p,\infty}(\R^d)}\, \left \lbrace \displaystyle \int_0^t \left (s^{\frac{1}{p}}\, f^\#(t) \right)^{ q_1} \,\displaystyle \frac{ds}{s}  \right \rbrace^{\frac{1}{q_2}}\\
\\
&= & \| f \|^{\frac{q_2-q_1}{q_2}}_{L^{p,\infty}(\R^d)}\, \| f \|^{\frac{ q_1}{q_2}}_{L^{p,q_1}(\R^d)}\\ \\
& \leq & \left (\displaystyle \frac{p}{q_1}\right)^{\frac{1}{q_1}}\,  \| f \|{L^{p,q_1}(\R^d)}\, \| f \|^{\frac{ q_1}{q_2}}_{L^{p,q_1}(\R^d)}\\
\end{array}
$$

\noindent since:

$$ \| f \|{L^{p,\infty}(\R^d)}=\displaystyle \sup_{t>0}\,t^{\frac{1}{p}}\, f^\#(t)$$

\end{proof}

\vskip 1cm

\begin{pte}

\noindent Let us consider $ 0 \leq p_1 < + \infty$, $ 0 \leq p_2 < + \infty$,  $ 0 \leq q_  < + \infty$. Then:

$$ p_1 \geq p_2 \Rightarrow L^{p_2,q }(\R^d) \hookrightarrow L^{p_1,q }(\R^d)$$


\end{pte}

\vskip 1cm

\begin{proof}\cite{BennettSharpley}$\,$\\

\noindent The secondary exponent is not involved, in so far as those inclusions are like the ones of the Lebesgue spaces $L^p$, and depend on the structure of the underlying measure space.

\end{proof}

\vskip 1cm

\section{New properties}
\begin{pte}(Cl. David)\\
\label{New}
\label{Max}
For any~$(f,g)$ belonging to~$L^{p,q}(\R^d) \times L^{p,q}(\R^d)$:

$$|||f+g|||_{L^{p,q}(\R^d)} \geq \max \, \left ( |||f |||_{L^{p,q}(\R^d)},|||g|||_{L^{p,q}(\R^d)} \right)$$

\end{pte}

\vskip 1cm

\begin{proof}

\noindent The proof relies on the following result:\\

\vskip 1cm

\begin{theorem}\label{Etage}
\noindent Let us denote by~$f$ a real-valued, positive function, defined on a measurable subset~$A$ of $\R^d$. Then, $f$ is the limit of an increasing sequence of positive simple function:

$$f= \displaystyle \lim_{n \to + \infty} \displaystyle\sum_{i\in I_n} f_{i,n}\,\mathbbm{1}_{A_{i,n}}$$

\noindent where, for any natural integer~$n$, $I_n$ is a countable set, and where, for any~$i$ belonging to~$I_n$, $f_{i,n}$ is a positive number, $A_{i,n}$ a measurable set, and $\mathbbm{1}_{A_{i,n}}$ is the characteristic function of the set~$A_{i,n}$.

\end{theorem}

\vskip 1cm

\begin{proof}

\noindent One requires just to examine the case of two simple functions, of the form:

$$f=   a_1 \,\1_{I_1  } \quad, \quad g=  a_2 \,\1_{I_2  } $$

\noindent where~$I_1 $ and~$I_2 $ are intervals of~$\R $, and~$a_1$, $a_2$, real numbers such that~$a_1< a_2$.\\

\noindent Then:\\

\begin{enumerate}

\item[$\rightsquigarrow$] If~$I_1$ and $I_2$ are disjoints, due to

$$  f^\# =  a_1\, \1_{[0,m_{1 }[}  $$

$$  g^\# =  a_2\, \1_{[0,m_{2 }[}  $$

\noindent one has:

$$ (f+g)^\#  \geq \max\, (f^\#,g^\#)$$

\item[$\rightsquigarrow$] If~$I_1=I_2$, one goes back, for~$f+g$, to a function of the form:

$$ \left  (\displaystyle \sum_{i=1}^2 a_i \right )\,\1_{I_1} $$

\noindent which leads to:

$$  f^\# =  a_1\, \1_{[m_{0},m_{1 }[}  $$

$$  g^\# =  a_2\, \1_{[m_{0},m_{1 }[}  $$

$$ (f+g)^\# =\left (\displaystyle \sum_{i=1}^2 a_i \right)\, \1_{[0 ,m_1[}  \geq \max\, (f^\#,g^\#) $$

\end{enumerate}

\begin{figure}[h!]
 \center{\psfig{height=8cm,width=7cm,angle=0,file=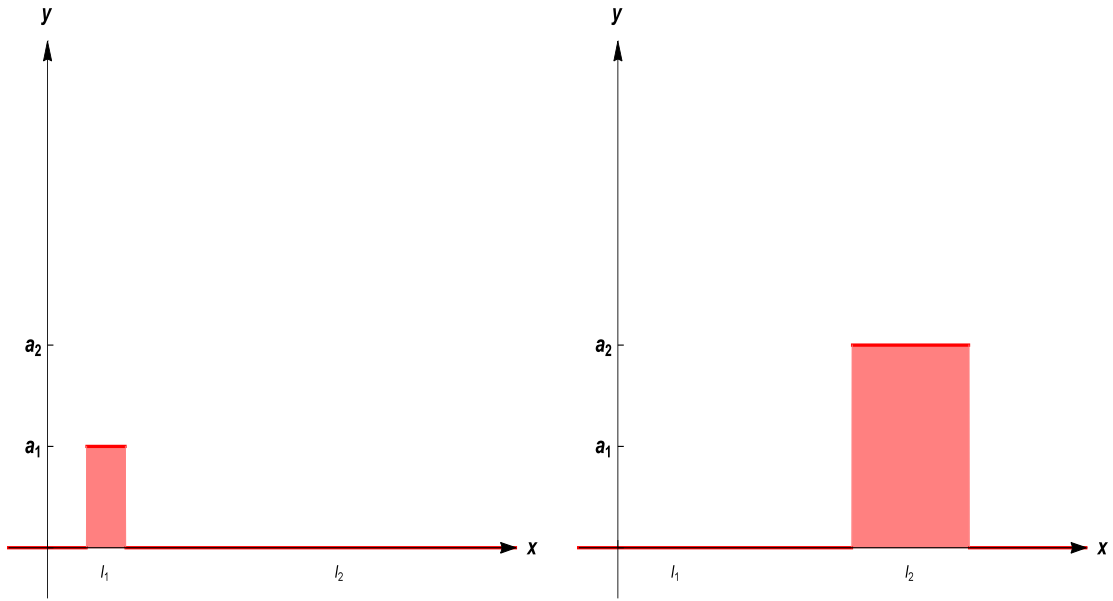}\hskip 0.5cm \psfig{height=8cm,width=7cm,angle=0,file=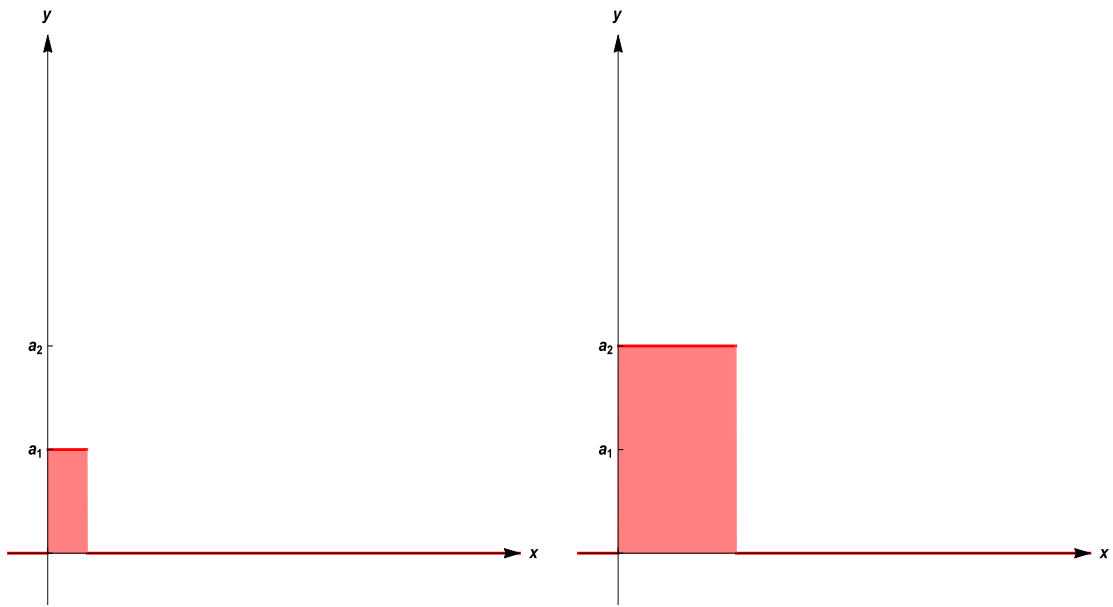}}
  \caption{The graph of $f $ (right), and $g $ (left).}
\end{figure}

\begin{figure}[h!]
 \center{\psfig{height=8cm,width=7cm,angle=0,file=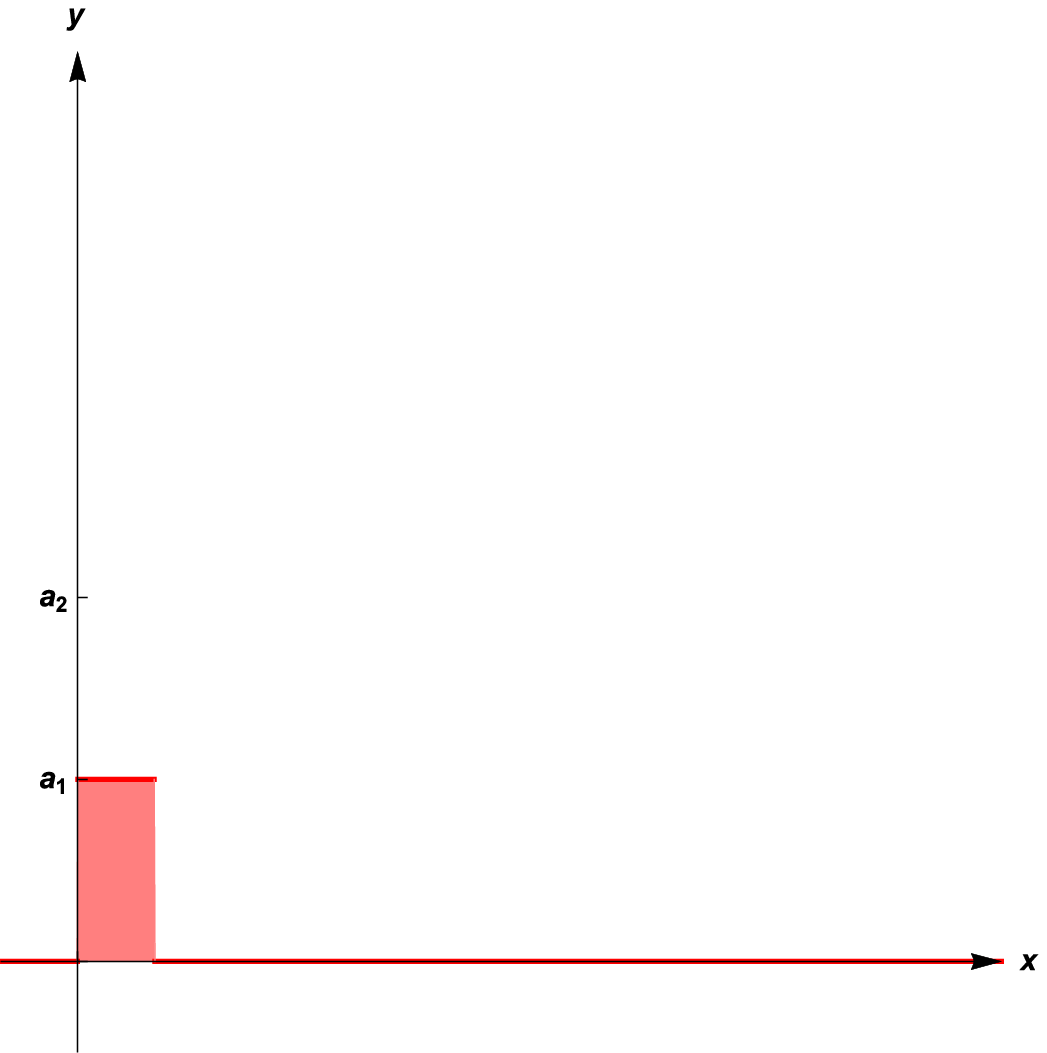}\hskip 0.5 cm \psfig{height=8cm,width=7cm,angle=0,file=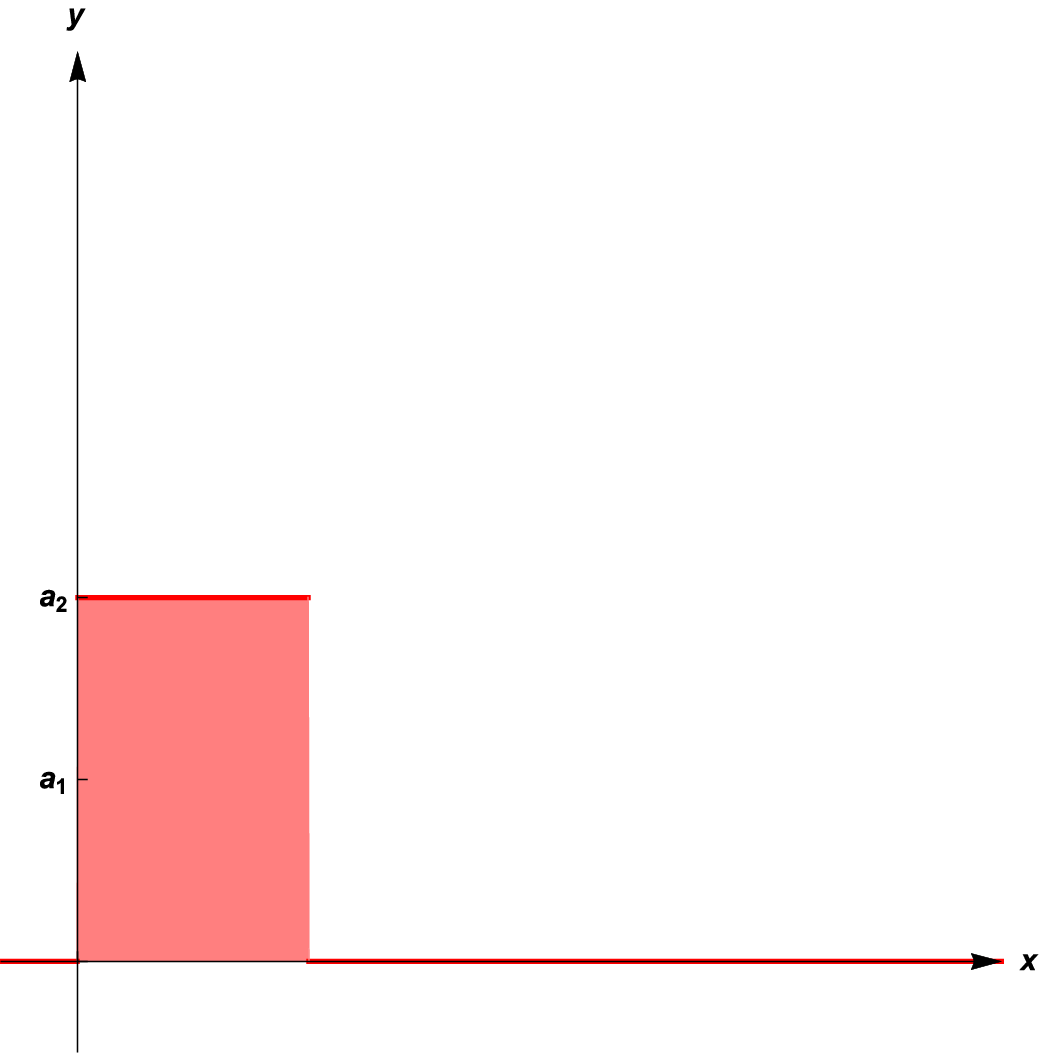}}
 \caption{The graph of $f^\#$ (right), and $g^\#$ (left).}
\end{figure}

\begin{figure}[h!]
 \center{\psfig{height=8cm,width=7cm,angle=0,file=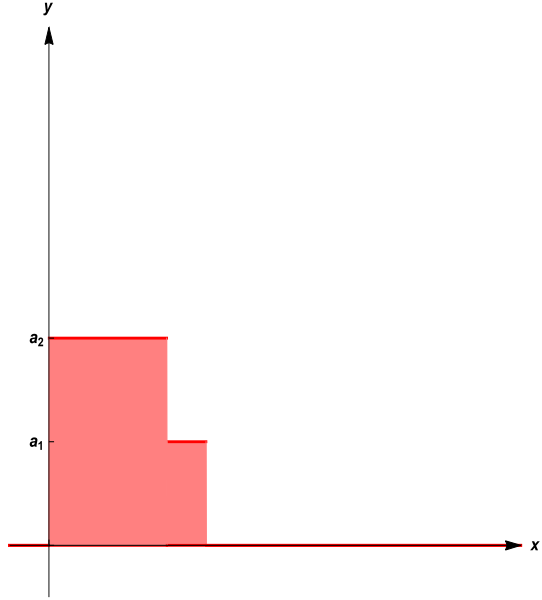}}
 \caption{The graph of $(f+g)^\#$.}
\end{figure}

\end{proof}








\vskip 1cm





\noindent Thus, for any strictly positive number~$t$:

$$(f+g)^\#  \geq \max\, (f^\#, g^\#)$$

\noindent and:

$$(f+g)^{\star \star}(t)=\displaystyle \frac{1}{t}\,\displaystyle \int_0^t (f+g)^\#(s)\,ds
\geq  \displaystyle \frac{1}{t}\,\displaystyle \int_0^t \max\, (f^\#(s), g^\#(s)) \,ds
  $$

\noindent Then:

$$(f+g)^{\star \star}(t) \geq  \max\, (f^{\star \star}(t), g^{\star \star}(t))$$

\noindent and:

$$t^{\frac{1}{p}}\, (f+g)^{\star \star}(t) \geq  \max\, t^{\frac{1}{p}}\, (f^{\star \star}(t), g^{\star \star}(t))$$

\noindent For~$q \geq 2$:

$$\left (t^{\frac{1}{p}}\, (f+g)^{\star \star}(t)\right)^q \geq  \max\, \left (\left (t^{\frac{1}{p}}\,  f^{\star \star}(t)\right)^q, \left (t^{\frac{1}{p}}\, g^{\star \star}(t)\right)^q \right)$$

\noindent which yields:

$$|||f+g|||_{L^{p,q}(\R^d)} \geq \max \, \left ( |||f |||_{L^{p,q}(\R^d)},|||g|||_{L^{p,q}(\R^d)} \right)$$

\end{proof}

\vskip 1cm

\vskip 1cm

\begin{lemma}\textbf{Maximal function related to a product of functions (Cl. David)}\\
\label{New2}

\noindent Let us denote by~$f$ and~$g$ two real-valued, measurable function, finite a.e. Then:

$$\forall \,t> 0 \, \quad \left (f\,g\right)^{\star \star} (t) \leq
\displaystyle\frac{2}{t}\,\displaystyle \int_0^{t}  f^\#(s) \,g^\#(s) \,ds
\leq
\displaystyle\frac{2}{t}\,\displaystyle \int_0^{t}  f^{\star \star} (s) \,g^{\star \star} (s) \,ds$$

\end{lemma}

\vskip 1cm

\begin{proof}
For any positive number~$s$:

$$\left (f\,g\right)^\# (s)  \leq  f^\#\left (\displaystyle \frac{s}{2}\right) \,g^\#\left (\displaystyle \frac{s}{2}\right)  $$

\noindent then, for any~$t>0$ :

$$\displaystyle\frac{1}{t}\,\displaystyle \int_0^t \left (f\,g\right)^\# (s)  \,ds \leq
\displaystyle\frac{1}{t}\,\displaystyle \int_0^t  f^\#\left (\displaystyle \frac{s}{2}\right) \,g^\#\left (\displaystyle \frac{s}{2}\right) \,ds
=
\displaystyle\frac{2}{t}\,\displaystyle \int_0^{\frac{t}{2}}  f^\#(s) \,g^\#(s) \,ds
\leq
\displaystyle\frac{2}{t}\,\displaystyle \int_0^{\frac{t}{2}}  f^{\star \star}(s) \,g^{\star \star}(s) \,ds$$

\noindent i.e.:

$$\left (f\,g\right)^{\star \star} (t) \leq
\displaystyle\frac{2}{t}\,\displaystyle \int_0^{t}  f^\#(s) \,g^\#(s) \,ds$$

\end{proof}

\vskip 1cm

\bibliographystyle{alpha}
\bibliography{BibliographieClaireRearrangement}

\end{document}